\documentclass[10pt,journal,compsoc]{IEEEtran}
\usepackage{booktabs} 
\usepackage{lipsum}   
\usepackage{ifpdf}  
\usepackage{amsmath}
\usepackage{graphicx}
\usepackage{threeparttable}
\usepackage{enumerate}
\usepackage{natbib}
\usepackage{url}  

\usepackage[colorlinks]{hyperref} 
\usepackage{float} 
\usepackage{verbatim} 
\usepackage{apalike} 
\restylefloat{figure} 
\restylefloat{table} 
\usepackage{geometry}
\usepackage{indentfirst}
\usepackage{color}
\usepackage{comment}
\usepackage{natbib}
\usepackage{bm}
\usepackage{rotating} 
\usepackage{longtable}
\usepackage{supertabular}
\usepackage{booktabs}
\usepackage{multirow}
\usepackage[table,xcdraw]{xcolor}
\usepackage{graphicx} 
\usepackage{algorithm}
\usepackage{algorithmic}

\usepackage{amsthm}
\usepackage{amsmath}

\usepackage{xcolor}  
\usepackage{caption} 
\usepackage{listings}
\usepackage{amssymb}
\usepackage{booktabs} 

\hypersetup{ 
    colorlinks=true,
    linkcolor=pptpurple,
    citecolor=pptpurple,
    urlcolor=pptpurple
}
\usepackage{pifont}
\definecolor{darkgreen}{RGB}{0, 176, 80}
\definecolor{islrorange}{RGB}{205, 100, 32}
\definecolor{islrblue}{RGB}{93, 165, 216}
\definecolor{whaleblue}{RGB}{51, 51, 178}
\definecolor{highlightblue}{RGB}{68, 68, 255}
\definecolor{pptorange}{RGB}{237, 125, 49}
\definecolor{pptgreen}{RGB}{112, 173, 71}
\definecolor{pptpurple}{RGB}{114, 51, 162}
\definecolor{highlightred}{RGB}{204, 0, 0}

\newtheorem{theorem}{\textbf{Theorem}}

\theoremstyle{plain}  

\theoremstyle{remark}

\newcommand{\cmark}{\ding{51}}%
\newcommand{\xmark}{\ding{55}}%

\begin{document}

\title{Parallel Network Reconstruction with Multi-directional Regularization}

\author{Zhaoyu Xing \IEEEmembership{}
and Wei Zhong \IEEEmembership{}

\IEEEcompsocitemizethanks{\IEEEcompsocthanksitem Zhaoyu Xing is the corresponding author with Department of Applied and Computational Mathematics and Statistics, University of Notre Dame, IN 46556, U.S.A.\protect \\
E-mail: zxing@nd.edu
\IEEEcompsocthanksitem Wei Zhong is with Department of Statistics \& Data Science, Xiamen University, Xiamen, Fujian 361005, China.}
}

%
%

\markboth{For arXiv}%
{Xing \MakeLowercase{\textit{et al.}}: Parallel Network Reconstruction with
Multi-directional Regularization}

\IEEEtitleabstractindextext{%
\begin{abstract}
Reconstructing large-scale latent networks from observed dynamics is crucial for understanding complex systems. However, the existing methods based on compressive sensing are often rendered infeasible in practice by prohibitive computational and memory costs. To address this challenge, we introduce a new distributed computing framework for efficient large-scale network reconstruction with parallel computing, namely PALMS (Parallel Adaptive Lasso with Multi-directional Signals). The core idea of PALMS is to decompose the complex global problem by partitioning network nodes, enabling the parallel estimation of sub-networks across multiple computing units. This strategy substantially reduces the computational complexity and storage requirements of classic methods. By using the adaptive multi-directional regularization on each computing unit, we also establish the consistency of PALMS estimator theoretically. Extensive simulation studies and empirical analyses on several large-scale real-world networks validate the computational efficiency and robust reconstruction accuracy of our approach.
\end{abstract}

}

\maketitle

\IEEEdisplaynontitleabstractindextext

%
\IEEEpeerreviewmaketitle


\section{Introduction}
\label{ch1}

Networks are pervasive in our world, characterizing interactions among diverse entities across numerous domains, such as the social networks with corresponding interactions between people and international trade networks with related economic growth. As network data and network-based dynamics have been widely collected \citep{chen2018structural}, the analysis of the complex dynamic systems with networks become one of the popular topics in both physics and machine learning \citep{HX2015, SD2021}. 

In many empirical conclusions, network structures are found profound to influence global dynamics and evolution of dynamic systems \citep{dhar2014prediction}, such as economic networks \citep{amiti2007economic, bramoulle2009identification}, social networks \citep{Santos2006, santos2008social, Choi2015}, and competitive interactions in evolutionary ultimatum games \citep{Antinyan2020}. In the classic network analysis, the network structures are commonly treated as fixed components in statistical model, such as the spatial auto-regressive models \citep{zou2021network, TZY2023} community detection methods \citep{Zhao2012, LeiRinaldo, airoldi2008mixed} and the models for identification  of network effect \citep{GG2010, Choi2015} and nodal effects \citep{xing2024golfs}. Consequently, all such analyses presuppose known network structures as prerequisites for model fitting, statistical inference and prediction.

However, the real-world large-scale network structures are usually latent and not accessible due to several reasons, such as budget constrains \citep{xing2025calms}, which makes all methods above inapplicable in practice. Instead of accessing network structures directly, networks reconstruction provides another alternative way to address this issue by inferring the latent networks utilizing the network-based dynamics \citep{HX2015}. As the observed network-based dynamics are usually noisy and affected by the latent structures, the network reconstruction are also known as a challenging ``inverse problem" \citep{HECKER200986, Timme_2014}.   

To address the network reconstruction problem, several approaches leveraging compressive sensing have recently been proposed \citep{Wang2011, HX2015, SD2021, adaptivesignallasso}. This group of methods reformulates the network reconstruction task as an optimization problem combining both the latent networks structure and the observed dynamics, and estimates the set of edges by solving the specially-designed optimization problem with respective to the adjacency matrix of the latent networks. Thus, different regularization can be used within this framework to obtain the sparse network estimations, such as the Lasso penalties \citep{HX2015} and Signal Lasso penalties \citep{SD2021}. Although several modifications based on compressive sensing are more accurate with appropriate theoretical properties \citep{xing2025calms}, they all require a complex transformation of the dynamic data and conduct the optimization algorithms based on the special data transformation, which can introduce a very large-scale loading matrix when the size of latent network is large  \citep{HX2015}. Due to the considerable algorithmic complexity and extensive memory requirements that scale with the size of the latent network, compressive-sensing based methods are largely infeasible for large-scale networks, such as online social networks and bio-networks. 

 
Motivated by these limitations while inferring the large-scale networks, we propose a new distributed framework to solve the network reconstruction problem, namely Parallel Adaptive Lasso with Multi-directional Signals (PALMS). PALMS decomposes the complex network reconstruction problem into many tasks via splitting of network dynamics, and reconstructs the sub-networks with parallel computing. We prove that the PALMS estimator using specially designed adaptive multi-directional penalization is consistent for latent true networks. Furthermore, we demonstrate that the proposed distributed framework is compatible with a wide range of network reconstruction methods that use compressive sensing. This allows for a significant reduction in computing time while maintaining a comparable level of accuracy. The comparisons between PALMS and previous classic methods in network reconstructions are shown in Table \ref{tab:method_comparison}.


\begin{table*}[b]
    \centering
    \caption{Comparison of Different Network Reconstruction Methods}
    \label{tab:method_comparison}
    \footnotesize
    \begin{tabular*}{\textwidth}{@{\extracolsep{\fill}}lccccc}
        \toprule
        \textbf{Method} & \textbf{Unbiasedness} & \textbf{Consistency} & \textbf{Fast Computing} & \textbf{Parallel Computing} & \textbf{Reference} \\
        \midrule
        Lasso        & \xmark & \xmark & \cmark & \xmark & \cite{HX2015} \\
        Signal Lasso & \xmark & \xmark & \xmark & \xmark & \cite{SD2021} \\
        ALMS         & \cmark & \cmark & \xmark & \xmark & \cite{xing2025calms} \\
        CALMS        & \cmark & \cmark & \xmark & \xmark & \cite{xing2025calms} \\
        PALMS (Ours) & \cmark & \cmark & \cmark & \cmark & This work \\
        \bottomrule
    \end{tabular*}
\end{table*}

The rest of the article is organized as follows. Section \ref{PRE} introduces the background of network reconstruction problem with two classic network dynamic models.  \ref{ModelSetting} presents the details of distributed network reconstruction algorithm and the theorem of PALMS estimator. The excellent finite-sample performance of PALMS is demonstrated based on simulations in Section \ref{Sim}. The real-data analysis with large-scale networks is shown in Section \ref{EA}. We conclude our contribution and make a discussion in Section \ref{Conclustion}. Technical conditions and the proof of theorem are given in the Appendix.

\section{Preliminaries}
\label{PRE}
To better introduce the methodology and algorithm of newly proposed PALMS, we introduce the background and preliminaries of network reconstructions and network-based dynamics in this section. 

In the field of network reconstruction, the sequences of nodal features are commonly treated as the dynamics generated with the fixed but latent networks.  To better describe the mechanism of network-based dynamics, several representative dynamic models for network-based systems are commonly considered in network reconstruction, such as evolutionary games \citep{Wang2011WNKY} and synchronization models \citep{Acebrón2005}. In many empirical analysis, these two models are widely proved to have the power to explain complex dynamic systems \citep{TGly2023, Wang2011WNKY, Frey, Nowak}. We next introduce the details of evolutionary games and synchronization models, and show that both of them can be included as special cases of our model proposed in Section \ref{ModelSetting}.  

Evolutionary game theory provides a framework for studying strategic interactions within populations of agents that adapt their behavior over time. While originally developed for economic studies, the evolutionary game is now widely used in fields such as biology and sociology \citep{Hummert2014}. The ultimatum game in networks is one of the well-known evolutionary game model for studying the dynamics of players' strategies and deposits within a network structure with exploration within the intricate fabric of large-scale networks, which are used to offer insights into how social structures shape economic and social behaviors and thus widely considered in network reconstruction problem \citep{HX2015, SD2021} as an representative dynamics. 

In the classic ultimatum game, participants are represented as nodes in a network, and the edges represent the interaction relationships between pairs of participants. In each round of the ultimatum game, one participant of a pair acts as the proposer and another as the responder, they need to decide the fraction to share a given amount of deposits without any communication. Specifically, the proposer will propose a fraction $p \in [0,1]$ for sharing, that is only $1-p$ of the deposit is shared to the responder. And then responder need to decides whether to accept or reject the offer based on a pre-set minimum accept fraction $q$. If the proposal is rejected, both of them gain $0$ from this round of games. A fair game process in considered in network-based ultimatum game, where both participants take turns playing the roles of proposer and responder for each round of the game. Then for each round of the balanced ultimatum game within one pair of participants, there are four potential outcomes summarized in the following \eqref{egePALMSSS} 
\begin{equation}
\label{egePALMSSS}
P_{i j}=\left\{\begin{array}{lll}
p_{j}+1-p_{i}, & p_{i} \geq r_{j}, & p_{j} \geq r_{i}; \\
1-p_{i}, & p_{i} \geq r_{j}, & p_{j}<r_{i}; \\
p_{j}, & p_{i}<r_{j}, & p_{j} \geq r_{i};\\
0, & p_{i}<r_{j}, & p_{j}<r_{i}.
\end{array}\right.
\end{equation}
where we have two participants $i$ and $j$ and their strategy pair denoted as $(p_i, r_i)$ and $(p_j, r_j)$, respectively.
As all participants are located in a network as nodes and each edge represents the balanced two-plays, one node (participant) with degree $d$ will play $d$-round of balanced ultimatum game independently with $d$ different direct neighbors, each of which contains two separate games as described above. Thus, the total income of the node (participant) in $t$-th round of game is the summation of the benefits from each edge, that is
\begin{equation}
\label{UGtotalY222P}
Y_i^t = \sum_{j=1}^k  a_{ij}P^t_{ij}, ~  i=1,\ldots,k,
\end{equation}
where $a_{ij}$ is the element of adjacency matrix, which represents the network structure of the ultimatum game. As shown in \ref{UGtotalY222P}, the network structure plays an important role in the network dynamics. The network reconstruction methods aim to estimate the latent network structures based only on the observed nodal dynamics. In the ultimatum game, latent large-scale networks refer to the systems comprising a substantial number of nodes in  societal, biological, and economic systems with large amount of people, genes and companies. The reconstruction of large-scale networks based on ultimatum game dynamics becomes an challenging topic in practice. 

Kuramoto model \citep{Kuramoto1975} is a physical model originally used to describes the synchronization of coupled oscillators, and is recently known for explaining the transition from incoherence to synchronized behavior in several complex systems, such as neuronal networks, power grids, and even flocks of birds \citep{hines2010topological,coutinho2013kuramoto}. Thus, it also widely considered as one of representative models in network reconstruction \citep{adaptivesignallasso, xing2025calms}. The model considers $N$ coupled oscillators, where each oscillator $i$ is characterized by its phase $\theta_i(t)$ and natural frequency $\omega_i$. The dynamics of the Kuramoto model are described by the following set of differential equations: 
\begin{equation}
\label{KM}
    \frac{d\theta_i}{dt} = \omega_i +   \sum_{j=1}^{N} a_{ij} \sin(\theta_j - \theta_i), \quad i = 1, 2, \dots, N,
\end{equation}
where $a_{ij}$ is from the adjacency matrix representing the coupling relationships between oscillators, and we assume the strength of interaction between the oscillators are identical. The right hand side of equation (\ref{KM}) is the rate of change for each nodes, which is effected by the direct neighbors of the oscillator $a_{ij}$ and the difference in their phase $\sin(\theta_j - \theta_i)$. When we treat the records of both the phase and it's rate of change as sequences of nodal features, the network structure will determine the nodal dynamics. As there are usually large amount of entities considered in the applications of Kuramoto model, the reconstruction of latent large-scale networks for complex dynamic system based on Kuramoto model is usually challenging.

In previous studies in economics and physics about evolutionary games and synchronization models, the networks are usually given and known as a part of experimental design. However, as we introduced in Section \ref{ch1}, the network structures that have significant effects in many real-world dynamic systems are commonly latent and network reconstruction approaches address the challenging inverse problem to infer the large-scale network structure based only on the observed nodal dynamics. In the next section, we formally introduce the problem of network reconstructions and propose our methodology of PALMS, where both ultimatum game \eqref{UGtotalY222P} and Kuramoto model \eqref{KM} introduced in this section are included as the special cases of our model. 



\section{Methodology and Algorithm} 
\label{ModelSetting}

To better introduce the methodology, we define the following notations for the problem setup of network reconstruction. Denote the adjacency matrix of the latent network $\mathcal{G}$ with node set $\mathcal{V}$ and edge set $\mathcal{E}$ as $\bm{A}\in \{0,1\}^{N\times N}$, and the size of networks is defined as the number of nodes included, that is $|\mathcal{V}|=N$.  Denote the observed nodal dynamic as $\bm{\Psi}\in\mathbb{R}^{N \times N}$ and nodal response as $Y \in \mathbb{R}^{N \times 1}$, and $\bm{\epsilon} \in \mathbb{R}^{N \times 1}$ represents the independent random errors with mean zero. Then the network dynamics for $t$-th round of network dynamics could be expressed as 
\begin{equation}
\label{tdm} 
    Y^t =  ( \bm{A} \circ \bm{\Psi}^t) \bm{1} + \bm{\epsilon}^t
\end{equation}
where $\circ$ is Hadamard product and $\bm{1}=[1,\cdots,1]'$. Note that  the evolutionary games in \eqref{UGtotalY222P} is a special case of (\ref{tdm}) with $\bm{\Psi}^t = \{P_{ij}^t\}$, and the Kuramoto model in \eqref{KM} is a special case of (\ref{tdm}) with $Y^t = [\Delta \theta^t/h]$ and $\bm{\Psi}^t = \{ \sin(\theta_i^t - \theta_j^t)\}$, $h$ is sampling interval. Both of the evolutionary games and the Kuramoto model are representative examples for network dynamics \citep{Wang2011, HX2015, SD2021,xing2025calms}. 

\begin{figure}[h]
    \centering
    \includegraphics[width=1\linewidth]{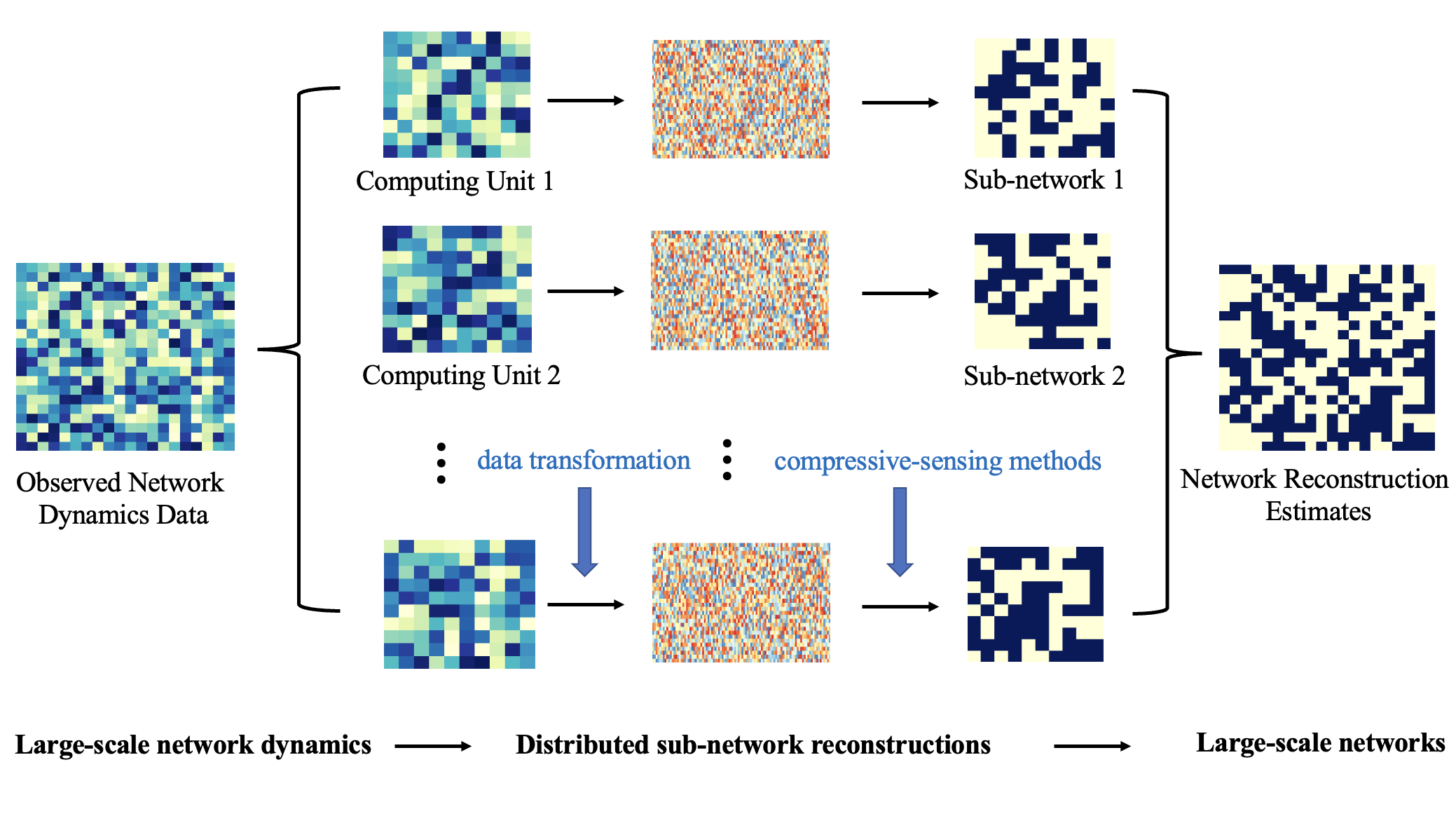}
    \caption{Computing process of distributed network reconstruction}
    \label{DistributedAlgorithmPlot}
\end{figure}

Drawing inspiration from the idea of cross-validation based on network node sampling \citep{NCV}, we propose a distributed computing framework for network reconstruction based on the random partitioning of nodes and corresponding network dynamics, as intuitively illustrated with the schematic diagram in Figure~\ref{DistributedAlgorithmPlot}. In the diagram, all matrix data are visualized using heatmaps. In practical applications, the observed network dynamics can consist of multiple matrices and one nodal dynamics is included, i.e., for each time $t$ in \eqref{tdm}, a corresponding network dynamic matrix $\bm{\Psi}^t$ is collected. After the node sampling, the dynamics of nodes in subnetworks can be transformed and the structures are estimated with compressive-sensing based methods on each computing units. By merging all the estimated structures, we have the whole network reconstructed.   

Specifically, we randomly allocate the nodes into $k$ groups, denoted as $\mathcal{G}_1, \cdots, \mathcal{G}_k$, where $\cup_{i=1}^k\mathcal{G}_i=\mathcal{V}$ and there is no overlap between different groups of nodes. According to the node partition, we divide the network dynamic data and transmit it to different computing units. Within each computing unit, a compressed sensing-based network reconstruction algorithm can be used to estimate the structure of the specific sub-network, which includes an initialization step involving data transformation and the optimization of an specially-designed objective function with regularization terms. Finally, by integrating the outputs from all computing units, we obtain a single reconstruction result for the entire large-scale network structure. Although the consistency of PALMS estimator are proved, we can uplift the accuracy and reduce the information loss with splitting by repeating the node partitioning in practice.

\begin{figure}[h]
    \centering
    \includegraphics[width=1\linewidth]{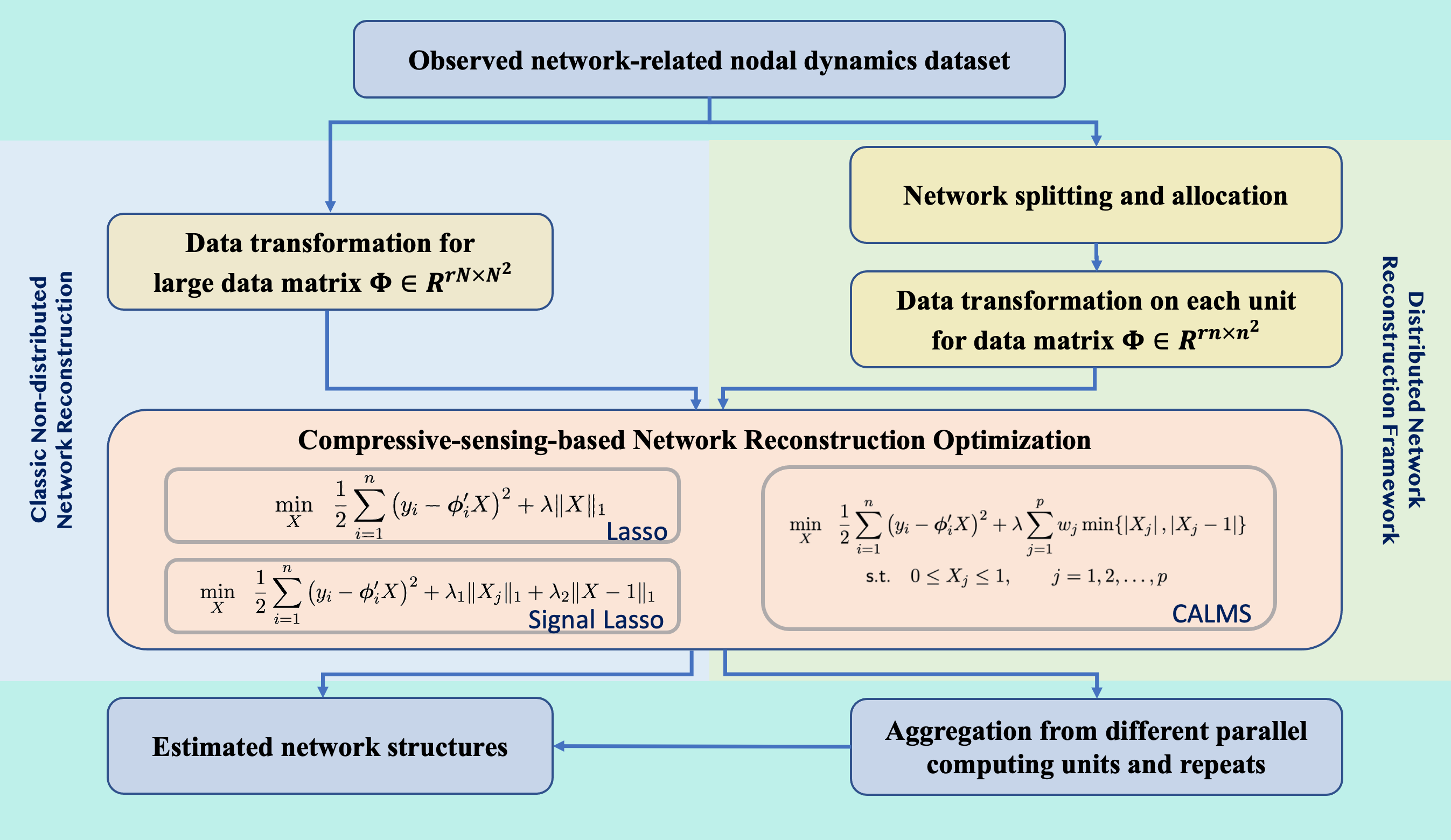}
    \caption{Comparison of the computing flows between distributed network reconstruction methods and classic  frameworks}
    \label{ComparisionOfFramework}
\end{figure}



For one network split, let $\mathcal{N}_{i},i\in[k]$ as a sequence of disjoint sets of nodes and $\bigcup_{i=1}^k \{ \mathcal{N}_{i}\}=\mathcal{V}$. Denoted the adjacency matrix for the sub-network between groups $\mathcal{N}_{k_1}$ and $\mathcal{N}_{k_2}$ as $\bm{A}^{\mathcal{N}_{k_1k_2}} \in \{0,1\}^{n\times n}$. The collected sub-network's dynamics with nodes $\mathcal{N}_{k_1k_2}$ consist the nodal response vector $Y^{\mathcal{N}_{k_1k_2}}$ and transformed design matrix $\bm{\Phi}^{\mathcal{N}_{k_1k_2}}$ of nodal dynamics (See \citealp{HX2015} for details of data-transformation). On each of the computing unit,  the subnetwork estimator based on compressed sensing can be expressed as the optimal solution to the following objective function
\begin{equation}
\label{overallest}
\tilde{\bm{A}}^{\mathcal{N}_{k_1k_2}} = \arg\min_{\bm{A}^{\mathcal{N}_{k_1k_2}}} \mathcal{L}(Y^{\mathcal{N}_{k_1k_2}}, \bm{\Phi}^{\mathcal{N}_{k_1k_2}}, \bm{A}^{\mathcal{N}_{k_1k_2}}),
\end{equation}
where $\mathcal{L}(\cdot)$ is a general loss function, typically including a loss metric related to residuals and a penalty term for network structures $\mathbf{A}$. Many existing network reconstruction methods \citep{Wang2011, HX2015, SD2021, adaptivesignallasso, xing2025calms} can be treated as special cases of Equation \eqref{overallest}. Denote the estimates of the entire network structure from $m$ rounds of node splitting as $\{\tilde{\bm{A}}_{(1)}, \cdots, \tilde{\bm{A}}_{(m)}\}$, then for any potential edges $A_{ij}$, the aggregate estimates is $\hat{A}_{ij} =   \sum_{l=1}^m [\tilde{A}_{ij}]_{(l)}/m.$ Correspondingly, the distributed estimator for the latent network can be written as:
\begin{equation}
\label{AGE}
\hat{\bm{A}} = \frac{1}{m} \sum_{l=1}^m \tilde{\bm{A}}_{(l)}.
\end{equation}
where $\hat{\bm{A}}_{(l)} = \bigcup_{k_1, k_2} \tilde{\bm{A}}^{\mathcal{N}_{k_1k_2}}_{(l)}$ is integrated of all subnetwork estimator in \ref{overallest}.


Note that many compressed sensing-based methods can be used to generate the subnetwork estimator $\tilde{\bm{A}}$ in \eqref{overallest}, and some special regularization possess desirable statistical properties, such as the consistency, for network reconstruction. 
By utilizing the new multi-directional regularization for subnetwork reconstruction, we proposed the distributed estimator for latent networks reconstruction as Parallel Adaptive Lasso with Multi-directional Signals (PALMS), which is obtained from
\begin{align}
&\tilde{\bm{A}}^{\mathcal{N}_{k_1k_2}}_{(l)} = \arg\min_{\bm{A}^{\mathcal{N}_{k_1k_2}}} \mathcal{L}^{(\text{PALMS})}(Y^{\mathcal{N}_{k_1k_2}} , \bm{\Phi}^{\mathcal{N}_{k_1k_2}}, \bm{A}^{\mathcal{N}_{k_1k_2}}),  \\
\label{PALMS_zhenghe}
& \hat{\bm{A}}_{(l)} = \bigcup_{k_1, k_2} \{\tilde{\bm{A}}^{\mathcal{N}_{k_1k_2}}_{(l)}\},\\
\label{PALMS}
& \hat{\bm{A}}^{(\text{PALMS})} = \frac{1}{m} \sum_{l=1}^m \hat{\bm{A}}_{(l)}.
\end{align}
where the subscript ${(l)}$ denotes the node dynamics and network structure corresponding to the $l$-th node partition, $\bigcup$ represents taking the union of the sub-network estimates to obtain the reconstruction result for the entire network, and the loss function $\mathcal{L}(\cdot)$ is
\begin{equation*}
\mathcal{L}^{(\text{PALMS})}(Y , \bm{\Phi}, \bm{A}) = \sum_{t} \| \bm{Y}^t - (\bm{A} \circ \bm{\Psi}^t) \|^2_F + P_\lambda(\bm{A})
\end{equation*}
where $P_\lambda(\bm{A})$ is the adaptive multi-directional penalty
\begin{equation}
\label{almsPen}
P_\lambda(\bm{A}) = \lambda_r \sum_{i,j} w_{ij} \min \{ |A_{ij}| , |A_{ij}-1| \},
\end{equation}
and $w_{ij}$ are adaptive weights for edge $A_{ij}$. A detailed discussion of the adaptive weights can be found in \citet{Zou2006} and \citet{xing2025calms}.

The following Theorem \ref{thm:consistency} demonstrates that the distributed network estimator $\hat{\bm{A}}_{PALMS}$ by is consistent for true network structures $\bm{A}^*$ when number of dynamics collected goes to infinity. By utilizing the new multi-directional regularization for subnetwork reconstruction on each computing unit, the PALMS method can consistently reconstruction the true latent networks while large reducing the computing time with parallel computing.

\begin{theorem} 
\label{thm:consistency}
Let $\bm{A}^*$ be the true adjacency matrix for the latent network. Under the Assumption (A1)-(A6), the aggregated estimator $\hat{\bm{A}}_{PALMS}$ defined in \eqref{PALMS} is a consistent estimator for $\bm{A}^*$ as $r \to \infty$:
$$
\hat{\bm{A}} \overset{p}{\to} \bm{A}^* \quad \text{as } r \to \infty.
$$
\end{theorem}


\begin{algorithm}
\caption{Distributed Network Reconstruction Algorithm}
\label{DistributedCALMS}
\begin{algorithmic}[1]

    \FOR{$l= 1, \cdots, m$}
        \STATE Randomly partition the network's node set $\mathcal{V}$ into $k$ groups, denoted as $\mathcal{G}_1, \cdots, \mathcal{G}_k$, such that $\cup_{i=1}^k\mathcal{G}_i=\mathcal{V}$;
        \FOR{$k_1, k_2 = 1, \cdots, k$ in parallel}
            \STATE Extract sub-networks $\mathcal{G}_{k_1}$ and $\mathcal{G}_{k_2}$ and their associated network dynamics; transmit to a sub-computational node;
            \STATE Denote the node response as $\bm{Y}^{\mathcal{N}_{k_1k_2}}$ and transform the node dynamics $\bm{\Psi}^{\mathcal{N}_{k_1k_2}}$ into the design matrix $\bm{\Phi}^{\mathcal{N}_{k_1k_2}}$ required by the compressed sensing-based algorithm;
            \STATE Estimate the sub-network structure by solving the optimization problem \eqref{overallest}, denoting the estimate as $\tilde{\bm{A}}^{\mathcal{N}_{k_1k_2}}_{(l)}$;
            \STATE Calculate the residuals $\mathbf{R}_{(l)} = \bm{Y}_{(l)}^{\mathcal{N}_{k_1k_2}} - \text{vec}(\tilde{\bm{A}}^{\mathcal{N}_{k_1k_2}}_{(l)})\bm{\Phi}^{\mathcal{N}_{k_1k_2}}$;
        \ENDFOR
        \STATE Obtain the integrated estimate for the entire network, $\hat{\bm{A}}_{(l)}$, according to \eqref{PALMS_zhenghe}.
    \ENDFOR
    \STATE Compute the final aggregated estimate using the formula in (\ref{PALMS}).
\end{algorithmic}
\end{algorithm}


We show the difference between the computational flow of the distributed network reconstruction framework and traditional methods visually in the flowchart in Figure \ref{ComparisionOfFramework}. For a large-scale network with $N$ nodes with $r$ rounds of network dynamics, the computational complexity of existing network reconstruction algorithms (e.g., Signal Lasso and CALMS) is $O(rN^3)$. Therefore the classic methods are almost not applicable for large-scale networks with long dynamic sequences. Furthermore, as detailed in \citet{HX2015}, these classical methods require storing the transformed design matrix $\bm{\Phi} \in \mathbb{R}^{rN \times N^2}$ during the initialization phase, which imposes a substantial memory requirement and makes it almost infeasible to solve large-scale network reconstruction problems on personal computers. 

As illustrated in Figure \ref{ComparisionOfFramework}, the proposed distributed computing framework introduces multi-step procedures for network splitting, data transmission, and parallel computation, which is applicable to most of network reconstruction algorithms based on compressed sensing techniques, such as the $L_1$-norm penalized method \citep{HX2015}, Signal Lasso \citep{SD2021} and CALMS \citep{xing2025calms}. Since only a small task reconstructing an $n \times n$ sub-network is allocated to each computing unit, the distributed framework can substantially reduce both the space and computational complexity of compressed sensing-based network reconstruction into $O(rn^3)$. We compare the practical performance of different network reconstruction methods in the following Section \ref{Sim} and show the uplifted computing speed of PALMS and the classic methods with proposed distributed framework.


\section{Simulation}    
\label{Sim}
In this section, we generate the network dynamics from classic and representative processes with different settings to test the effectiveness of proposed distributed framework for network reconstruction. 


To evaluate and compare the performance of different methods in network reconstruction, we generate the true network structures based on the Erdős-Rényi (ER) random networks model \citep{Erdos1959} which is one of the representative random network models and commonly used in the simulation of network reconstruction problem \citep{SD2021, adaptivesignallasso, xing2025calms}. In the ER random network model, each pair of nodes is connected by an edge with a pre-set probability $p$. This results in a binomial degree distribution where the degree of each node follows the binomial distribution. 


To fully simulate the real-world complex dynamic systems, we refer to previous works on network reconstruction \citep{HX2015, SD2021, xing2025calms} and generate $r$ rounds of network-based dynamics from the three classic models under the given network structures: the multi-normal dynamics, the evolutionary ultimatum game dynamics and the synchronization model dynamics. The background of the models is introduced in Section \ref{PRE} and details of the generating process are described as below.

\begin{table*}[b]
    \begin{threeparttable}
    \centering
    \caption{The comparisons based on ER random networks and different types of network dynamics}
    \label{GDPtable}
    \footnotesize
    \begin{tabular*}{\textwidth}{@{\extracolsep{\fill}}clcccc}
        \toprule
        \textbf{DGP} & \multicolumn{1}{c}{\textbf{Methods}} & \textbf{CPU Time} & \textbf{MSE} & \textbf{SRNL} & \textbf{SREL} \\
        \midrule
        \multirow{6}{*}{DGP 1}
         & Lasso          & 0.127 (0.055)      & 0.501 (0.014)      & 0.996 (0.008) & 0.000 (0.000) \\
         & P-Lasso        & 0.400 (0.081)      & 1.550 (0.396)      & 0.936 (0.013) & 0.000 (0.000) \\
         & Signal Lasso   & 82.373 (6.282)     & 0.375 (0.022)      & 0.385 (0.236) & 0.273 (0.212) \\
         & P-Signal Lasso & 26.673 (2.917)     & 0.473 (0.012)      & 0.808 (0.025) & 0.243 (0.022) \\
         & CALMS          & 1089.263 (113.364) & 0.435 (0.013)      & 0.561 (0.018) & 0.568 (0.015) \\
         & PALMS          & 43.779 (6.355)     & 0.462 (0.011)      & 0.539 (0.017) & 0.536 (0.016) \\
        \midrule
        \multirow{6}{*}{DGP 2}
         & Lasso          & 0.120 (0.026)      & 0.247 (0.012)      & 0.998 (0.004) & 0.000 (0.000) \\
         & P-Lasso        & 0.401 (0.077)      & 0.675 (0.182)      & 0.950 (0.010) & 0.000 (0.000) \\
         & Signal Lasso   & 77.153 (6.202)     & 0.226 (0.016)      & 0.719 (0.097) & 0.069 (0.018) \\
         & P-Signal Lasso & 26.720 (1.263)     & 0.327 (0.019)      & 0.799 (0.031) & 0.285 (0.042) \\
         & CALMS          & 1096.677 (121.47)  & 0.433 (0.02)       & 0.571 (0.020) & 0.554 (0.030) \\
         & PALMS          & 33.676 (5.166)     & 0.467 (0.011)      & 0.508 (0.014) & 0.608 (0.030) \\
        \midrule
        \multirow{6}{*}{DGP 3}
         & Lasso          & 0.141 (0.041)      & 0.729 (0.015)      & 1.000 (0.001) & 0.000 (0.000) \\
         & P-Lasso        & 0.271 (0.085)      & 298.289 (814.566)  & 0.965 (0.014) & 0.000 (0.000) \\
         & Signal Lasso   & 51.807 (6.851)     & 0.365 (0.39)       & 0.244 (0.032) & 0.766 (0.087) \\
         & P-Signal Lasso & 10.652 (2.091)     & 0.289 (0.015)      & 0.100 (0.032) & 0.938 (0.021) \\
         & CALMS          & 568.436 (113.454)  & 0.361 (0.022)      & 0.359 (0.028) & 0.742 (0.024) \\
         & PALMS          & 14.130 (2.760)     & 0.271 (0.015)      & 0.004 (0.003) & 0.999 (0.001) \\
        \bottomrule
    \end{tabular*}
    \begin{tablenotes}
        \item[\textit{Note:}] The parameters are set as $n=50$, $k=5$, $r=5$, the expectation of network density is $0.5$.
    \end{tablenotes}
    \end{threeparttable}
\end{table*}

\noindent\textbf{Data generating process 1}: The loading matrix $\bm{\Phi} $ follows a multivariate normal distribution $\mathcal{N}(\mu,\Sigma)$ where $\mu_j\sim \mathrm{Unif}(-1,1)$, $\Sigma_{jj} \sim \mathrm{Unif}(1,3)$ and $\Sigma_{ij}=0$ for $i\neq j$. We generate the response vector $Y^t=\{Y_1^t,\cdots,Y_2^t\}$ at time $t$ via \eqref{tdm} with $\epsilon_i \sim \mathcal{N}(0,1)$, where $t=1,\cdots,r$.

\noindent\textbf{Data generating process 2}: The data are generated based on the evolutionary ultimatum game in networks. The each player's strategy pair $(p_i, q_i)$ is initially generated randomly from a standard uniform distribution and updated based on their payoffs. Specifically, the probability that player $i$ adopts a new strategy is given by $\{1 + \exp(Y_i^r - Y_j^r/N)\}^{-1}$, where $Y_i^r$ and $Y_j^r$ denote the payoffs of players $i$ and $j$ in round $r$, and $N$ is the size of network. The payoff matrix $\bm{P}^{(r)}$ for the $r$-th round is computed based on Equation \eqref{egePALMSSS}. Then transfer the nodal dynamics into $\bm{\Phi}_i $ matrix with same procedure as in \citet{HX2015}, whose $r$-th row corresponds to the payoff coefficients from the payoff matrix $P^{(r)}$ in round $r$, and obtain the nodal response vectors in $r$-round of dynamics from \eqref{tdm} as $Y_i = \left( Y_i^1, Y_i^2, \ldots, Y_i^r \right)^\top$ with $\epsilon_i$ from standard normal distribution $\mathcal{N}(0,1)$. That is, the payoffs of player $i$ is obtained from 
\[
\bm{Y} = \bm{\Phi} X + \bm{\Xi},
\]
where: $Y$ is a stacked vector of all players' payoffs, represented as $
Y = (Y_1', Y_2', \dots, Y_N')' \in \mathbb{R}^{rN \times 1}$,
where $Y_1, Y_2, \dots, Y_N$ are individual payoffs of the players, $X = vec(\bm{A})$. The matrix $\bm{\Phi}$ is a block-diagonal matrix composed of the individual payoff coefficient matrices, expressed as $\bm{\Phi}  = \mathrm{diag}(\bm{P}^{(1)}, \bm{P}^{(2)}, \dots, \bm{P}^{(r)}) \in \mathbb{R}^{rN \times N(N-1)}$. Finally, $\bm{\Xi}$ is the vector of random error terms. This formulation captures how each player's payoffs are influenced by their connections with other players and the strategies employed over multiple rounds.

\noindent\textbf{Data generating process 3}: The data are generated based on a synchronization model introduced in Section \ref{PRE}. We consider the discrete case and fix time step $h=0.01$, and generate the data with the following update formula for the $t$-th step:
\begin{equation}
\frac{y_i^t}{h}= c \sum_{j=1}^{N} a_{i j}  \sin \left(\theta_{j}^t-\theta_{i}^t \right),
\end{equation}
where $c$ is a constant for coupling strength. For $r$ discrete changes of $N$ vertices, we define $\phi_{i j}^t=c \sin \left(\theta_{j}^t-\theta_{i}^t\right)$ and let $Y_i=(y_i^1,\ldots,y_i^r)$ and $Y=(Y_1',\ldots,Y_N')$, $X_i = (A_{i1},\ldots,A_{iN})$ and $X=(X_1',\ldots,X_N')$,
\begin{equation}
\bm{\Phi}_{i}=\left(\begin{array}{cccc}
\phi_{i 1}^1 & \phi_{i, 2}^1 & \cdots & \phi_{i_{N}}^1 \\
\phi_{i 1}^2 & \phi_{i, 2}^2 & \cdots & \phi_{i_{N}}^2 \\
\vdots & \vdots & \vdots & \vdots \\
\phi_{i 1}^r & \phi_{i, 2}^r & \cdots & \phi_{i_{N}}^r
\end{array}\right)
\end{equation}
and $\bm{\Phi}=\mathrm{diag}(\bm{\Phi}_{1},...,\bm{\Phi}_{N})$. The final simulated data are generated from $Y = \bm{\Phi}  X  + \bm{\Xi}$.

\begin{table*}[b]
    \begin{threeparttable}
    \centering
    \caption{The comparisons based on ER random networks and different sizes of networks}
    \label{Ntable}
    \footnotesize 
    \begin{tabular*}{\textwidth}{@{\extracolsep{\fill}}llcccc}
        \toprule
        \textbf{Size} & \textbf{Method} & \textbf{CPU Time} & \textbf{MSE} & \textbf{SRNL} & \textbf{SREL} \\
        \midrule
        \multirow{6}{*}{30}
         & Lasso          & 6.112 (10.354)      & 1.214 (1.241)     & 0.997 (0.009) & 0.000 (0.000) \\
         & P-lasso        & 0.706 (0.390)       & 0.720 (0.024)     & 0.964 (0.020) & 0.040 (0.012) \\
         & Signal Lasso   & 20.052 (8.362)      & 0.268 (0.018)     & 0.138 (0.036) & 0.663 (0.083) \\
         & P-Signal Lasso & 9.495 (3.968)       & 0.302 (0.028)     & 0.149 (0.042) & 0.891 (0.040) \\
         & CALMS          & 114.274 (62.561)    & 0.376 (0.014)     & 0.376 (0.033) & 0.711 (0.017) \\
         & PALMS          & 23.753 (11.185)     & 0.317 (0.015)     & 0.187 (0.030) & 0.857 (0.016) \\
        \midrule
        \multirow{6}{*}{40}
         & Lasso          & 19.906 (142.343)    & 0.752 (0.044)     & 0.999 (0.002) & 0.000 (0.000) \\
         & P-lasso        & 0.778 (0.395)       & 0.713 (0.016)     & 0.976 (0.013) & 0.027 (0.005) \\
         & Signal Lasso   & 48.197 (21.491)     & 0.286 (0.014)     & 0.108 (0.015) & 0.593 (0.043) \\
         & P-Signal Lasso & 18.009 (7.041)      & 0.304 (0.017)     & 0.121 (0.021) & 0.912 (0.019) \\
         & CALMS          & 600.076 (325.875)   & 0.391 (0.018)     & 0.380 (0.029) & 0.695 (0.024) \\
         & PALMS          & 45.683 (21.068)     & 0.324 (0.014)     & 0.178 (0.021) & 0.863 (0.015) \\
        \midrule
        \multirow{6}{*}{50}
         & Lasso          & 120.437 (952.318)   & 0.886 (0.313)     & 0.996 (0.004) & 0.000 (0.000) \\
         & P-lasso        & 0.950 (0.515)       & 0.722 (0.007)     & 0.981 (0.005) & 0.026 (0.005) \\
         & Signal Lasso   & 152.129 (8.081)     & 0.298 (0.009)     & 0.096 (0.014) & 0.572 (0.027) \\
         & P-Signal Lasso & 33.254 (14.683)     & 0.300 (0.011)     & 0.119 (0.018) & 0.908 (0.018) \\
         & CALMS          & 2907.002 (1546.948) & 0.391 (0.011)     & 0.358 (0.023) & 0.699 (0.012) \\
         & PALMS          & 104.204 (50.769)    & 0.332 (0.011)     & 0.177 (0.028) & 0.844 (0.020) \\
        \bottomrule
    \end{tabular*}
    \begin{tablenotes}
        \item[\textit{Note:}] The parameters are set as $k=5$ and $r=10$, the expected network density is $0.5$. The network dynamics are generated by evolutionary game.
    \end{tablenotes}
    \end{threeparttable}
\end{table*}

To compare the contributions of distributed computation in network reconstruction, we refer to the classic work in network reconstruction and compare the performance of different methods in the simulations with following 4 measures: (1) CPU Time: CPU time used for completing network reconstruction estimations by each algorithm in seconds; (2) Mean Square Error (MSE): This metric is used to measure the magnitude of the overall network reconstruction error, calculated using the following formula $\text{MSE}(\hat{X})=\sum_{j=1}^p(X_j-\hat{X}_j)^2/p$; (3) Success Rates for Nonexistent Links (SRNL): This metric is used to measure the accuracy of the overall network reconstruction for empty edges, calculated by the ratio between the number of correctly predicted nonexistent links and the total number of nonexistent links: $\text{SRNL}(\hat{X})= \sum_{j=1}^p \mathrm{I}(X_j=0) \mathrm{I}(\hat{X}_j=0)/\sum_{j=1}^p \mathrm{I}(X_j=0)$; 
(4) Success Rates for Existing Links (SREL): this metric is used to measure the accuracy of the overall network reconstruction for existing edges, calculated by the ratio between the number of correctly predicted links and the total number of links $\text{SREL}(\hat{X})= \sum_{j=1}^p \mathrm{I}(X_j=1)  \mathrm{I}(\hat{X}_j=1)/\sum_{j=1}^p \mathrm{I}(X_j=1).$

For each simulation settings, we repeats for $500$ times and report the averages and corresponding standard deviations of $4$ measures in tables. The details of settings and comparisons are shown in the footnotes of each table. As the distributed algorithm in Algorithm \ref{DistributedCALMS} aims mainly on addressing the computing challenges of network reconstruction, we focus on the uplift of computing speed to achieve similar accuracy of reconstructions. We compare the performance between classic network reconstruction method Lasso \citep{HX2015}, Signal Lasso \citep{SD2021}, CALMS \citep{xing2025calms} and their performance with parallel computing with Algorithm \ref{DistributedCALMS}, which are named as ``P-Lasso", ``P-Signal Lasso " and ``PALMS", respectively.

The comparison of reconstruction performance across different network dynamic models using ER random networks is given in Table \ref{GDPtable}. The distributed PALMS algorithm demonstrates significantly shorter computation times than CALMS across all types of network dynamics. While achieving high computational efficiency, PALMS maintains reconstruction accuracy comparable to CALMS in terms of Mean Squared Error (MSE). Although Lasso method can be computed efficiently with well-developed R packages \citep{JSSv033i01}, it has no power to estimate the edges exactly with SREL$=0$ across all cases. 

The comparison of the finite-sample performance of various network reconstruction methods across different network scales using ER random networks are summarized in Table \ref{Ntable}. Distributed framework demonstrates substantial advantages in computation time for most of the network reconstruction methods. As network size increases, traditional reconstruction methods like Signal Lasso and CALMS exhibit the largely increasing computation times, which indicates their computational sensitivity to the size of latent networks. In contrast, computational advantage of the distributed framework becomes more pronounced for larger networks, which is consistent with our analysis of complexity in Section \ref{ModelSetting}. 

 \begin{table*}[b]
    \begin{threeparttable}
    \centering
    \caption{The performance of network reconstruction with increasing rounds of network dynamics}
    \label{RoundTable}
    \footnotesize
    \begin{tabular*}{\textwidth}{@{\extracolsep{\fill}}cccccc}
        \toprule
        \textbf{Method} & \textbf{r} & \textbf{CPU Time} & \textbf{MSE}  & \textbf{SRNL} & \textbf{SREL} \\
        \midrule
        \multirow{4}{*}{ALMS}  & 3  & 43.705 (7.064)   & 0.297 (0.011) & 0.800 (0.012) & 0.281 (0.033) \\
                               & 5  & 74.392 (10.696)  & 0.275 (0.018) & 0.820 (0.018) & 0.308 (0.033) \\
                               & 10 & 172.677 (29.999) & 0.255 (0.014) & 0.833 (0.011) & 0.365 (0.04)  \\
                               & 20 & 360.514 (65.108) & 0.219 (0.012) & 0.860 (0.009) & 0.435 (0.030) \\
        \midrule
        \multirow{4}{*}{PALMS} & 3  & 4.942 (0.671)    & 0.56 (0.026)  & 0.347 (0.034) & 0.843 (0.031) \\
                               & 5  & 5.666 (1.046)    & 0.55 (0.028)  & 0.355 (0.037) & 0.863 (0.036) \\
                               & 10 & 5.754 (0.904)    & 0.437 (0.019) & 0.534 (0.024) & 0.690 (0.035) \\
                               & 20 & 6.848 (2.064)    & 0.374 (0.017) & 0.577 (0.022) & 0.839 (0.025) \\
        \bottomrule
    \end{tabular*}
    \begin{tablenotes}
        \item[\textit{Note:}] The parameters are set as $n=50$, $k=5$, the expectation of network density is $0.1$.
    \end{tablenotes}
    \end{threeparttable}
\end{table*}

As shown in Table \ref{RoundTable}, both CALMS and PALMS exhibit decreasing MSE as the number of network dynamic rounds $r$ increases. This empirically validates the consistency results established in Theorem \ref{thm:consistency} that when large amount of network-based dynamics are collected, we can consistently reconstruct the latent true networks with our proposed methods. Meanwhile, as shown in Table \ref{RoundTable}, the computing time of distributed method PALMS are relatively not sensitive to long sequences of dynamics. Together with the empirical results in Table \ref{GDPtable}, the distributed framework in \eqref{PALMS} can handle the challenging task of reconstruction of large-scale networks with long sequences. In the next section, we demonstrate the usefulness of distributed network with several large-scale real-world datasets where the non-distributed approaches are not feasible.

\section{Empirical analysis}
\label{EA}

\begin{figure}[ht]
    \centering
    \includegraphics[width=1\linewidth]{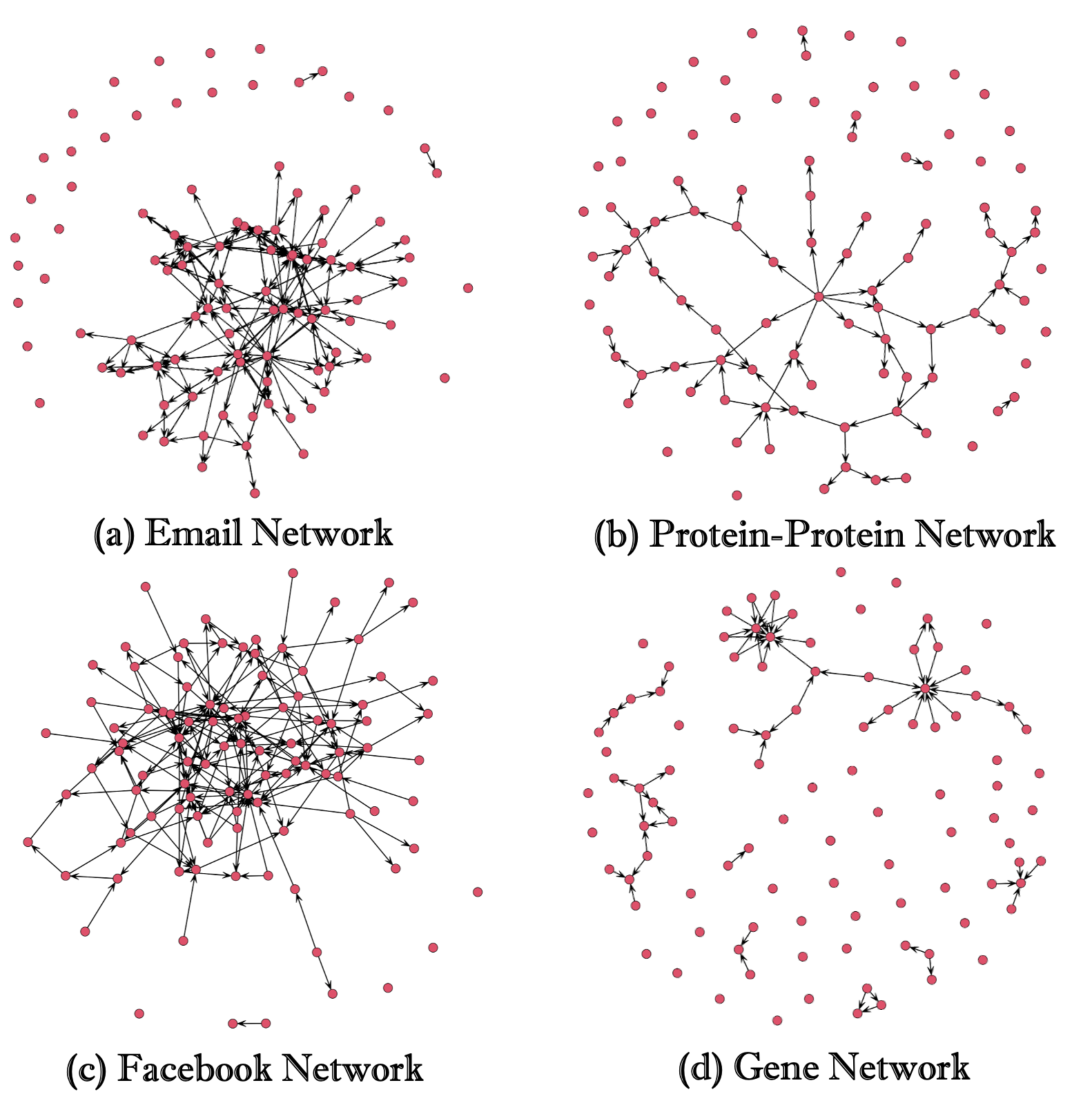}
   \caption{The typical network structures of real-world datasets}
   \label{EANETS}
\end{figure}

To better show the effectiveness of our proposed distributed framework, we use classic methods combining with the distributed framework in Algorithm \ref{DistributedCALMS} and PALMS to infer four real-world structures: the email-connection network among the scholars in European institutions\footnote{https://snap.stanford.edu/data/email-Eu-core.html}, the gene interaction network in biological experiments\footnote{https://snap.stanford.edu/data/S-cerevisiae.html}, 
the Protein-protein Interaction Network\footnote{http://vlado.fmf.uni-lj.si/pub/networks/ data/bio/Yeast/Yeast.htm}
and the Facebook friendship network\footnote{https://snap.stanford.edu/data/ego-Facebook.html}. 
As Figure \ref{EANETS} shows, four real-world datasets have significant difference in network structures. More details of the real-world network data are given in Section \ref{DEA} of Appendix. With real structures of networks are obtained from the dataset, we generate $10$ rounds of network-based dynamics from the evolutionary game introduced in Section \ref{Sim} for network reconstruction task. As the classic non-distributed network reconstruction methods are time-consuming and almost not feasible for these networks, we only compare the performance in network reconstruction with different methods combined with proposed distributed framework based on the real-world networks. The tuning process and evaluations are the same as demonstrated in Section \ref{Sim}.

\begin{table*}[b]
    \centering
    \caption{Empirical analysis results on real-world networks}
    \label{Eempiricalresulttable}
    \footnotesize
    \begin{tabular*}{\textwidth}{@{\extracolsep{\fill}}lcccc}
        \toprule
        \textbf{Dataset} & \textbf{Method} & \textbf{SREL} & \textbf{SRNL} & \textbf{MSE} \\
        \midrule
        \multirow{4}{*}{Email network}
         & Lasso          & 0.0000 & 1.0000 & 0.1105 \\
         & P-lasso        & 0.0000 & 0.9051 & 7.1571 \\
         & P-Signal Lasso & 0.5260 & 0.5831 & 0.4232 \\
         & PALMS          & 0.9377 & 0.0638 & 0.8396 \\
        \midrule
        \multirow{4}{*}{Protein-protein interaction network}
         & Lasso          & 0.0000 & 1.0000 & 0.0080 \\
         & P-lasso        & 0.0000 & 1.0000 & 0.0080 \\
         & P-Signal Lasso & 0.5685 & 0.7222 & 0.2790 \\
         & PALMS          & 0.7109 & 0.5747 & 0.4242 \\
        \midrule
        \multirow{4}{*}{Facebook network}
         & Lasso          & 0.0000 & 1.0000 & 0.0207 \\
         & P-lasso        & 0.0000 & 0.9993 & 0.0225 \\
         & P-Signal Lasso & 0.7231 & 0.5808 & 0.4163 \\
         & PALMS          & 0.8925 & 0.3535 & 0.6353 \\
        \midrule
        \multirow{4}{*}{Gene network}
         & Lasso          & 0.0000 & 1.0000 & 0.0078 \\
         & P-lasso        & 0.0000 & 1.0000 & 0.0078 \\
         & P-Signal Lasso & 0.7488 & 0.8398 & 0.1609 \\
         & PALMS          & 0.8293 & 0.7752 & 0.2243 \\
        \bottomrule
    \end{tabular*}
\end{table*}

The estimation results of different distributed network reconstruction methods are given in Table \ref{Eempiricalresulttable}. For the reconstruction problem of the large-scale and sparse network in our empirical analysis, PALMS achieves the dominant performance in detecting the exist connections in large-scale network reconstructions in term of the SREL. Although MSE of Lasso methods are very small, SREL gives a straightforward measure of the ability to estimate the latent connections in large-scale and sparse networks, as the MSE can also be extremely small when an empty network is estimated for sparse networks with few edges. 


As Figure \ref{EANETS} and Table \ref{SEAD} shown, the datasets collected are large-scale with complex structures of edges. While several previously proposed methods are no longer feasible for inferring the latent network structures, the proposed distributed algorithm can reduce the computing complexity and requirements of memory, which can still handle the reconstruction task for large-scale network in practice with accurate estimation. 

\section{Conclusion and discussion} 
\label{Conclustion}
This paper introduced a distributed computing framework to effectively overcome the computational barriers in large-scale network reconstruction. By partitioning the network and employing parallel processing, PALMS makes computationally-intensive compressive sensing techniques feasible for large networks by drastically reducing execution time. We theoretically proved that PALMS estimator enjoys the appropriate statistical property consistency for the latent true network structure. The superior performance of PALMS was demonstrated in both extensive simulations and real-world data analysis, and the R tools for PALMS are publicly available online\footnote{https://github.com/ZhaoyuXingStat/PALMS}. Although our method achieves high accuracy with significantly less computation time, future work could investigate the impact of different network partitioning strategies and the design to a robust network reconstruction methods. 



\section*{Acknowledgments}
This research is supported by the National Natural Science Foundation of China  (11922117, 12231011, 71988101), the National Key R\&D Program of China (2022YFA1003802), and the National Statistical Science Research Grants of China (2022LD08).

\newpage
\small
\renewcommand{\theequation}{A.\arabic{equation}}
\setcounter{equation}{0}
\renewcommand{\thesubsection}{A.\arabic{subsection}}
\renewcommand{\thetheorem}{A.\arabic{theorem}}
\setcounter{theorem}{0}
  
 \setlength{\bibsep}{4pt plus 0.3ex} 
\bibliographystyle{chicago} 
\footnotesize
\bibliography{sample} 

\begin{thebibliography}{}

\bibitem[\protect\citeauthoryear{Acebr{\'o}n, Bonilla, P{\'e}rez~Vicente, Ritort, and Spigler}{Acebr{\'o}n et~al.}{2005}]{Acebrón2005}
Acebr{\'o}n, J.~A., L.~L. Bonilla, C.~J. P{\'e}rez~Vicente, F.~Ritort, and R.~Spigler (2005).
\newblock The kuramoto model: A simple paradigm for synchronization phenomena.
\newblock {\em Reviews of Modern Physics\/}~{\em 77\/}(1), 137--185.

\bibitem[\protect\citeauthoryear{Airoldi, Blei, Fienberg, and Xing}{Airoldi et~al.}{2008}]{airoldi2008mixed}
Airoldi, E.~M., D.~Blei, S.~Fienberg, and E.~Xing (2008).
\newblock Mixed membership stochastic blockmodels.
\newblock {\em Advances in Neural Information Processing Systems\/}~{\em 21}.

\bibitem[\protect\citeauthoryear{Amiti and Cameron}{Amiti and Cameron}{2007}]{amiti2007economic}
Amiti, M. and L.~Cameron (2007).
\newblock Economic geography and wages.
\newblock {\em The Review of Economics and Statistics\/}~{\em 89\/}(1), 15--29.

\bibitem[\protect\citeauthoryear{Antinyan, Horvath, and Jia}{Antinyan et~al.}{2020}]{Antinyan2020}
Antinyan, A., G.~Horvath, and M.~Jia (2020).
\newblock Positional concerns and social network structure: An experiment.
\newblock {\em European Economic Review\/}~{\em 129}, 103547.

\bibitem[\protect\citeauthoryear{Bramoull{\'e}, Djebbari, and Fortin}{Bramoull{\'e} et~al.}{2009}]{bramoulle2009identification}
Bramoull{\'e}, Y., H.~Djebbari, and B.~Fortin (2009).
\newblock Identification of peer effects through social networks.
\newblock {\em Journal of econometrics\/}~{\em 150\/}(1), 41--55.

\bibitem[\protect\citeauthoryear{Bu, Zhao, Cai, Xue, Zhu, Lu, Zhang, Sun, Ling, Zhang, et~al.}{Bu et~al.}{2003}]{bu2003topological}
Bu, D., Y.~Zhao, L.~Cai, H.~Xue, X.~Zhu, H.~Lu, J.~Zhang, S.~Sun, L.~Ling, N.~Zhang, et~al. (2003).
\newblock Topological structure analysis of the protein--protein interaction network in budding yeast.
\newblock {\em Nucleic Acids Research\/}~{\em 31\/}(9), 2443--2450.

\bibitem[\protect\citeauthoryear{Chen and Lei}{Chen and Lei}{2018}]{NCV}
Chen, K. and J.~Lei (2018).
\newblock Network cross-validation for determining the number of communities in network data.
\newblock {\em Journal of the American Statistical Association\/}~{\em 113\/}(521), 241--251.

\bibitem[\protect\citeauthoryear{Chen, An, An, Guan, and Hao}{Chen et~al.}{2018}]{chen2018structural}
Chen, Z., H.~An, F.~An, Q.~Guan, and X.~Hao (2018).
\newblock Structural risk evaluation of global gas trade by a network-based dynamics simulation model.
\newblock {\em Energy\/}~{\em 159}, 457--471.

\bibitem[\protect\citeauthoryear{Coutinho, Goltsev, Dorogovtsev, and Mendes}{Coutinho et~al.}{2013}]{coutinho2013kuramoto}
Coutinho, B., A.~Goltsev, S.~Dorogovtsev, and J.~Mendes (2013).
\newblock Kuramoto model with frequency-degree correlations on complex networks.
\newblock {\em Physical Review E—Statistical, Nonlinear, and Soft Matter Physics\/}~{\em 87\/}(3), 032106.

\bibitem[\protect\citeauthoryear{Dhar, Geva, Oestreicher-Singer, and Sundararajan}{Dhar et~al.}{2014}]{dhar2014prediction}
Dhar, V., T.~Geva, G.~Oestreicher-Singer, and A.~Sundararajan (2014).
\newblock Prediction in economic networks.
\newblock {\em Information Systems Research\/}~{\em 25\/}(2), 264--284.

\bibitem[\protect\citeauthoryear{Erd{\H{o}}s and R{\'e}nyi}{Erd{\H{o}}s and R{\'e}nyi}{1959}]{Erdos1959}
Erd{\H{o}}s, P. and A.~R{\'e}nyi (1959).
\newblock On random graphs.
\newblock {\em Publicationes Mathematicae\/}~{\em 6}, 290--297.

\bibitem[\protect\citeauthoryear{Frey}{Frey}{2010}]{Frey}
Frey, E. (2010).
\newblock Evolutionary game theory: Theoretical concepts and applications to microbial communities.
\newblock {\em Physica A: Statistical Mechanics and its Applications\/}~{\em 389\/}(20), 4265--4298.

\bibitem[\protect\citeauthoryear{Friedman, Hastie, and Tibshirani}{Friedman et~al.}{2010}]{JSSv033i01}
Friedman, J.~H., T.~Hastie, and R.~Tibshirani (2010).
\newblock Regularization paths for generalized linear models via coordinate descent.
\newblock {\em Journal of Statistical Software\/}~{\em 33\/}(1), 1–22.

\bibitem[\protect\citeauthoryear{Gallo and Syngjoo}{Gallo and Syngjoo}{2016}]{Choi2015}
Gallo, E. and C.~Syngjoo (2016).
\newblock Networks in the laboratory.

\bibitem[\protect\citeauthoryear{Ghiglino and Goyal}{Ghiglino and Goyal}{2010}]{GG2010}
Ghiglino, C. and S.~Goyal (2010).
\newblock Keeping up with the neighbors: social interaction in a market economy.
\newblock {\em Journal of the European Economic Association\/}~{\em 8\/}(1), 90--119.

\bibitem[\protect\citeauthoryear{Han, Shen, Wang, and Di}{Han et~al.}{2015}]{HX2015}
Han, X., Z.~Shen, W.-X. Wang, and Z.~Di (2015).
\newblock Robust reconstruction of complex networks from sparse data.
\newblock {\em Physical Review Letters\/}~{\em 114\/}(2), 028701.

\bibitem[\protect\citeauthoryear{Hecker, Lambeck, Toepfer, {van Someren}, and Guthke}{Hecker et~al.}{2009}]{HECKER200986}
Hecker, M., S.~Lambeck, S.~Toepfer, E.~{van Someren}, and R.~Guthke (2009).
\newblock Gene regulatory network inference: Data integration in dynamic models—a review.
\newblock {\em Biosystems\/}~{\em 96\/}(1), 86--103.

\bibitem[\protect\citeauthoryear{Hines, Blumsack, Sanchez, and Barrows}{Hines et~al.}{2010}]{hines2010topological}
Hines, P., S.~Blumsack, E.~C. Sanchez, and C.~Barrows (2010).
\newblock The topological and electrical structure of power grids.
\newblock In {\em 2010 43rd Hawaii International Conference on System Sciences}, pp.\  1--10. IEEE.

\bibitem[\protect\citeauthoryear{Hummert, Bohl, Basanta, Deutsch, Werner, Theißen, Schroeter, and Schuster}{Hummert et~al.}{2014}]{Hummert2014}
Hummert, S., K.~Bohl, D.~Basanta, A.~Deutsch, S.~Werner, G.~Theißen, A.~Schroeter, and S.~Schuster (2014).
\newblock Evolutionary game theory: cells as players.
\newblock {\em Molecular BioSystems\/}~{\em 10\/}(12), 3044–3065.

\bibitem[\protect\citeauthoryear{Kuramoto}{Kuramoto}{1975}]{Kuramoto1975}
Kuramoto, Y. (1975).
\newblock Self-entrainment of a population of coupled non-linear oscillators.
\newblock In {\em International Symposium on Mathematical Problems in Theoretical Physics}, pp.\  420--422.

\bibitem[\protect\citeauthoryear{Lei and Rinaldo}{Lei and Rinaldo}{2015}]{LeiRinaldo}
Lei, J. and A.~Rinaldo (2015).
\newblock Consistency of spectral clustering in stochastic block models.
\newblock {\em The Annals of Statistics\/}, 215--237.

\bibitem[\protect\citeauthoryear{Leskovec and Mcauley}{Leskovec and Mcauley}{2012}]{leskovec2012learning}
Leskovec, J. and J.~Mcauley (2012).
\newblock Learning to discover social circles in ego networks.
\newblock {\em Advances in Neural Information Processing Systems\/}~{\em 25}.

\bibitem[\protect\citeauthoryear{Nowak}{Nowak}{2006}]{Nowak}
Nowak, M.~A. (2006).
\newblock {\em Evolutionary dynamics: exploring the equations of life}.
\newblock Harvard University Press.

\bibitem[\protect\citeauthoryear{Santos, Pacheco, and Lenaerts}{Santos et~al.}{2006}]{Santos2006}
Santos, F.~C., J.~M. Pacheco, and T.~Lenaerts (2006).
\newblock Evolutionary dynamics of social dilemmas in structured heterogeneous populations.
\newblock {\em Proceedings of the National Academy of Sciences\/}~{\em 103\/}(9), 3490--3494.

\bibitem[\protect\citeauthoryear{Santos, Santos, and Pacheco}{Santos et~al.}{2008}]{santos2008social}
Santos, F.~C., M.~D. Santos, and J.~M. Pacheco (2008).
\newblock Social diversity promotes the emergence of cooperation in public goods games.
\newblock {\em Nature\/}~{\em 454\/}(7201), 213--216.

\bibitem[\protect\citeauthoryear{Shi, Hu, Jin, Shen, Tan, and Yu}{Shi et~al.}{2023}]{adaptivesignallasso}
Shi, L., J.~Hu, L.~Jin, C.~Shen, H.~Tan, and D.~Yu (2023).
\newblock Robust and efficient network reconstruction in complex system via adaptive signal lasso.
\newblock {\em Physical Review Research\/}~{\em 5\/}(4), 043200.

\bibitem[\protect\citeauthoryear{Shi, Shen, Jin, Shi, Wang, and Boccaletti}{Shi et~al.}{2021}]{SD2021}
Shi, L., C.~Shen, L.~Jin, Q.~Shi, Z.~Wang, and S.~Boccaletti (2021).
\newblock Inferring network structures via signal lasso.
\newblock {\em Physical Review Research\/}~{\em 3\/}(4), 043210.

\bibitem[\protect\citeauthoryear{Tho, Ding, Hui, Welsh, and Zou}{Tho et~al.}{2023}]{TZY2023}
Tho, Z.~Y., D.~Ding, F.~K. Hui, A.~Welsh, and T.~Zou (2023).
\newblock On the robust estimation of spatial autoregressive models.
\newblock {\em Econometrics and Statistics\/}, 2452--3062.

\bibitem[\protect\citeauthoryear{Timme and Casadiego}{Timme and Casadiego}{2014}]{Timme_2014}
Timme, M. and J.~Casadiego (2014, aug).
\newblock Revealing networks from dynamics: an introduction.
\newblock {\em Journal of Physics A: Mathematical and Theoretical\/}~{\em 47\/}(34), 343001.

\bibitem[\protect\citeauthoryear{Traulsen and Glynatsi}{Traulsen and Glynatsi}{2023}]{TGly2023}
Traulsen, A. and N.~E. Glynatsi (2023).
\newblock The future of theoretical evolutionary game theory.
\newblock {\em Philosophical Transactions of the Royal Society B\/}~{\em 378\/}(1876), 20210508.

\bibitem[\protect\citeauthoryear{Wang, Lai, Grebogi, and Ye}{Wang et~al.}{2011}]{Wang2011}
Wang, W.-X., Y.-C. Lai, C.~Grebogi, and J.~Ye (2011).
\newblock Network reconstruction based on evolutionary-game data via compressive sensing.
\newblock {\em Physical Review X\/}~{\em 1\/}(2), 021021.

\bibitem[\protect\citeauthoryear{Wang, Ni, Lai, and Grebogi}{Wang et~al.}{2011}]{Wang2011WNKY}
Wang, W.-X., X.~Ni, Y.-C. Lai, and C.~Grebogi (2011).
\newblock Pattern formation, synchronization, and outbreak of biodiversity in cyclically competing games.
\newblock {\em Physical Review E: Statistical, Nonlinear, and Soft Matter Physics\/}~{\em 83\/}(1), 011917.

\bibitem[\protect\citeauthoryear{Xing, Tan, Zhong, and Shi}{Xing et~al.}{2025}]{xing2025calms}
Xing, Z., H.~Tan, W.~Zhong, and L.~Shi (2025).
\newblock Calms: Constrained adaptive lasso with multi-directional signals for latent networks reconstruction.
\newblock {\em Neurocomputing\/}, 129545.

\bibitem[\protect\citeauthoryear{Xing, Wan, Wen, and Zhong}{Xing et~al.}{2024}]{xing2024golfs}
Xing, Z., Y.~Wan, J.~Wen, and W.~Zhong (2024).
\newblock Golfs: feature selection via combining both global and local information for high dimensional clustering.
\newblock {\em Computational Statistics\/}~{\em 39\/}(5), 2651--2675.

\bibitem[\protect\citeauthoryear{Zhao, Levina, and Zhu}{Zhao et~al.}{2012}]{Zhao2012}
Zhao, Y., E.~Levina, and J.~Zhu (2012).
\newblock Consistency of community detection in networks under degree-corrected stochastic block models.
\newblock {\em The Annals of Statistics\/}~{\em 40\/}(4), 2266--2292.

\bibitem[\protect\citeauthoryear{Zou}{Zou}{2006}]{Zou2006}
Zou, H. (2006).
\newblock The adaptive lasso and its oracle properties.
\newblock {\em Journal of the American Statistical Association\/}~{\em 101\/}(476), 1418--1429.

\bibitem[\protect\citeauthoryear{Zou, Luo, Lan, and Tsai}{Zou et~al.}{2021}]{zou2021network}
Zou, T., R.~Luo, W.~Lan, and C.-L. Tsai (2021).
\newblock Network influence analysis.
\newblock {\em Statistica Sinica\/}~{\em 31\/}(4), 1727--1748.

\end{thebibliography}

\newpage

\section*{Proof of Theorem \ref{thm:consistency}}
\label{PT}

\section*{Details of real datasets}
\label{DEA}


\begin{table*}[b]
    \centering
    \caption{Summary statistics of real-world networks used in empirical analysis}
    \label{SEAD}
    \footnotesize
    \begin{tabular*}{\textwidth}{@{\extracolsep{\fill}}lccc}
        \toprule
        \textbf{Real Dataset} & \textbf{Number of nodes} & \textbf{Number of edges} & \textbf{Number of triangles} \\
        \midrule
        Email-connection network & 1005 & 25571 & 105461 \\
        Gene interaction network & 690 & 1094 & 72 \\
        Protein-protein network & 2361 & 7182 & 536 \\
        Facebook friendship network & 4039 & 88234 & 1612010 \\
        \bottomrule
    \end{tabular*}
\end{table*}

The mail-networks dataset is a network dataset generated from email communications within a large European research institution. It provides anonymized information on email exchanges between the institution's members and highlights how people are connected through their email communications. This large-scale network have 1005 nodes and 25571 edges and 105461 triangles, which includes both internal and external communication links. 

The gene network dataset represents a network of gene regulation in the yeast Saccharomyces cerevisiae. In this network, each node corresponds to an operon, which is a group of genes transcribed together into a single mRNA molecule. The network includes 690 nodes, 1094 directed edges and 72 triangles. A directed edge from operon $i$ to operon $j$ means that $i$ is regulated by a transcriptional factor encoded by $j$.

The protein-protein interaction network in budding yeast dataset \citep{bu2003topological}, consists of a network that represents interactions between proteins in Saccharomyces cerevisiae (budding yeast). The dataset includes 2361 vertices (proteins) and 7182 edges (interactions), with 536 loops. 

Facebook friendship network is a collection of anonymity Facebook data that focuses on users' friends lists \citep{leskovec2012learning}. It contains ego networks, which capture the relationships between a user (ego) and their friends (alters), as well as the relationships between those friends. The large network are obtained by merging the ego networks and includes 4039 nodes and 88234 edges. The network has a high average clustering coefficient of 0.6055, indicating that users' friends tend to be friends with each other.

\end{document}